\documentclass[a4paper,11pt]{amsproc}
\usepackage{amssymb}
\usepackage{amsmath}
\usepackage{amsthm,lmodern} 
\usepackage{tikzsymbols}
\usepackage[foot]{amsaddr}

\usepackage{framed, booktabs, float}
\usepackage{xcolor}
\usepackage[backref]{hyperref}
\hypersetup{
    colorlinks,
    linkcolor={red!60!black},
    citecolor={green!60!black},
    urlcolor={blue!60!black}
}
\usepackage{enumerate}
\usepackage{mathtools}
\usepackage{url}
\usepackage[T1]{fontenc}
\usepackage{geometry}
\usepackage{tikz}
\usepackage{pgfplots}
\usetikzlibrary{3d}
\usetikzlibrary{calc}

\theoremstyle{plain}
\newtheorem{Theorem} {Theorem} [section]
\newtheorem{Proposition} [Theorem] {Proposition}
\newtheorem{Lemma} [Theorem] {Lemma}
\newtheorem{Corollary} [Theorem] {Corollary}

\newtheorem{Remark} [Theorem] {Remark}

\theoremstyle{definition}
\newtheorem{Definition} [Theorem] {Definition}

\newcommand{\Ff}{{\mathbb F}}
\newcommand{\F}{{\mathbb F}}

\newcommand{\cF}{{\mathcal F}}

\newcommand{\cP}{{\mathcal P}}

\newcommand{\eps}{\varepsilon}

\newcommand{\<}{\langle}
\renewcommand{\>}{\rangle} 
\renewcommand{\phi}{\varphi} 

\newcommand{\onum}{\sigma}

\DeclareMathOperator\p{\mathcal{P}}
\DeclareMathOperator\pg{\mathrm{PG}}
\DeclareMathOperator\h{\mathsf{H}}
\DeclareMathOperator\w{\mathsf{W}}
\DeclareMathOperator\q{\mathsf{Q}}

\title{New bounds and constructions for large partial $m$-ovoids and related structures}

\author[Bamberg]{John Bamberg$^1$}
\address{$^1$Department of Mathematics \& Statistics, 
The University of Western Australia, Perth, Australia.}

\author[Bishnoi]{Anurag Bishnoi$^2$}
\address{$^2$Delft Institute of Applied Mathematics, Delft University of Technology, Netherlands.}

\author[Ihringer]{Ferdinand Ihringer$^3$}
\address{$^3$Department of Mathematics, 
Southern University of Science and Technology,
Shenzhen, China.}

\email{john.bamberg@uwa.edu.au}
\email{A.Bishnoi@tudelft.nl}
\email{ferdinand.ihringer@gmail.com}
\email{A.ravi@tudelft.nl}

\author[Ravi]{Ananthakrishnan Ravi $^4$}
\address{$^4$Delft Institute of Applied Mathematics, Delft University of Technology, Netherlands.}

\begin{document}

\begin{abstract}
  We use $p$-rank bounds on partial ovoids and the classical bounds on Ramsey numbers to obtain upper bounds on the size of partial $m$-ovoids in finite classical polar spaces. These bounds imply a uniform non-existence result of $m$-ovoids over all families of finite classical polar spaces.
  
In the special case of the symplectic spaces over the binary field, we prove an equivalence between partial $m$-ovoids and a generalisation of Oddtown families from extremal set theory that has been studied under the name of $m$-nearly orthogonal sets. We give a new construction for large partial $2$-ovoids in these spaces and thus $2$-nearly orthogonal sets over the binary field. 
This construction uses triangle-free graphs associated to certain BCH codes whose complements have low $2$-rank and it gives an asymptotic improvement over the previous best constructions.
We give another construction of triangle-free graphs using a binary projective cap, which has low complementary rank over the reals.
This improves the bounds in the recently introduced rank-Ramsey problem of Beniamini, Linial, and Shraibman.
It also gives better constructions of large partial $m$-ovoids for $m > 2$ in the binary symplectic space.
\end{abstract}

\maketitle

\section{Introduction}
One of the main areas of research in finite geometry is to understand extremal substructures inside a given geometrical space. 
For example, understanding the largest possible sets of points in finite projective or affine spaces with the property that no three points are collinear is a classical problem with connections to additive combinatorics \cite{Alon13, EG17} and coding theory \cite{Hill78}. 
We will study problems of a similar flavour in the geometrical spaces known as the \textit{finite classical polar spaces}, which were introduced by Tits \cite{Tits1959} as natural geometries corresponding to the finite classical simple groups. 
Their main geometric property, as shown by Buekenhout and Shult \cite{BS74}, is the so-called \emph{one or all axiom}: if $P$ is a point and $\ell$ a line not incident with $P$, then 
$P$ is collinear with either all points of $\ell$ or with exactly one point of $\ell$. 
In fact, abstract polar spaces are defined as point-line geometries satisfying this axiom. 
We will be interested in finite classical polar spaces that arise from
equipping a vector space over a finite field with a sesquilinear or quadratic form. 
The points of the polar space are then the one-dimensional totally isotropic or singular vector subspaces and the lines are the two-dimensional totally isotropic or totally singular vector subspaces, with respect to the given form. 

Besides applications to group theory, polar spaces give rise to various interesting combinatorial objects like strongly regular graphs \cite{BvM}, two-weight codes \cite{CK86}, and Tur\'{a}n graphs \cite{Ball-Pepe}. 
The incidence structures known as generalised quadrangles are special cases of polar spaces, and since their inception, they have been used to construct various extremal objects in graph theory.
Recently, other polar spaces have also found combinatorial applications, with the symplectic space over the binary field being implicitly used in the work of Conlon and Ferber on multicolour Ramsey numbers \cite{CF21}, and graphs related to the rank one Hermitian spaces being used in the breakthrough work of Mattheus and Verstraete on off-diagonal Ramsey numbers \cite{MV23}. 

The study of ovoids in polar spaces dates back to the work of Tits on geometric constructions of certain Suzuki groups \cite{ThasHvM,Tits62,Tits62Suzuki}, by using a particular polar space. 
In \cite{Thas81}, Thas synthesised these objects as ovoids of polar spaces: a set $\mathcal{O}$ of points such that every maximal totally isotropic subspace meets $\mathcal{O}$ in exactly one point. 
In terms of the collinearity graph of the polar space,\footnote{The collinearity graph of any point-line geometry is the graph on points with two points adjacent if they lie on a common line.} these are independent sets that meet every maximal clique in exactly one vertex, and in particular, they are the largest possible independent sets. 
The study of $m$-ovoids begins implicitly with the work of Segre \cite{Segre65}, but was introduced for finite generalised quadrangles by Thas \cite{Thas:Interesting}, and further extended by Shult and Thas \cite{ST1996} to all finite polar spaces. 
Here, an $m$-ovoid is defined to be a set $\mathcal{O}$ of points such that every maximal totally isotropic subspace meets $\mathcal{O}$ in exactly $m$ points. 
A partial $m$-ovoid is a set of points that meets every maximal totally isotropic subspace in at most $m$ points. 
These objects are natural to study even in polar spaces that have no $m$-ovoids, since they measure how large a set can be while still avoiding $m+1$ points on a generator.

\textit{In this paper}, we show applications in the reverse direction to the applications mentioned above: we apply bounds and constructions from extremal combinatorics to improve our knowledge of partial $m$-ovoids in finite classical polar spaces. 
Our first contribution is a general Ramsey-theoretic upper bound on partial $m$-ovoids. 
It shows that every large partial $m$-ovoid gives a graph with no large clique and no large coclique, and hence that classical bounds on Ramsey numbers imply upper bounds on partial $m$-ovoids.

\begin{Theorem}\label{Ramseypartialmovoid}
Let $k$ be the maximum size of a partial $1$-ovoid in a finite polar space $\mathcal{P}$. 
Then for any partial $m$-ovoid $\mathcal{O}_m$ of $\mathcal{P}$, we have 
\[
|\mathcal{O}_m| < R(m+1,k+1).
\]
\end{Theorem}

Together with known $p$-rank bounds on partial ovoids, this gives general upper bounds for partial $m$-ovoids in all finite classical polar spaces. 
From these bounds we can deduce the following new non-existence result for $m$-ovoids.

\begin{Corollary}\label{Ramseypartialmovoid_corollary}
    An $m$-ovoid in a finite classical polar space of rank $r$ over a finite field of characteristic $p$ does not exist if 
    \[
m \log_p(2r+2+p) < \frac{r-1}{p-1}.
\]
\end{Corollary}
Although this bound on $m$ for the non-existence of $m$-ovoids is weaker than the earlier bounds in the cases where those bounds apply, it has the advantage of being applicable to \textit{all} finite classical polar spaces (see \cite{BKLP2007}).
Therefore, for the cases not covered by the previous results, where the earlier techniques cannot show any non-existence results, we get new non-existence results.
For comparison, Bamberg et al. \cite{BKLP2007} gave the following lower bound on $m$ for which an $m$-ovoid can exist in elliptic quadrics, symplectic spaces in even dimension, and Hermitian spaces in odd dimension:
\[
m\ge \frac{\sqrt{4 q^{e+r-1}+9}-3}{2 (q-1)}.
\]
Here, $q=p^h$ is the size of the defining field, where $p$ is a prime, and $r$ is the rank of the polar space.
The parameter $e\in\{1,\frac{3}{2},2\}$ depends on the type of the polar space; see Section \ref{background}. 
Asymptotically, we have $m=\Omega(q^{r/2})$ if we fix $q$ and vary $r$, and $m=\Omega\left(p^{h(e+r-3)/2}\right)$ if we fix $(r,h)$ and vary $p$. 
This was recently improved by De Beule, Mannaert, and Smaldore \cite{DeBMS}: they proved $m=\Omega( \sqrt{r} q^{r/2})$ if we fix $q$ and vary $r$, with no asymptotic change in the lower bound on $m$ if we fix $(r,h)$ and vary $p$. 
The results of \cite{BKLP2007}, \cite{DeBMS}, and \cite{HamiltonMathon} on the non-existence of $m$-ovoids only apply to the polar spaces $\q^-(2r+1,q)$, $\w(2r-1,q)$, and $\h(2r,q^2)$. 
In contrast, the Ramsey bound applies to all finite classical polar spaces, thus proving new non-existence results for the remaining infinite families of finite classical polar spaces (see Table~\ref{table:overview}).

In Section \ref{sec:spectralbounds}, we give other upper bounds on partial $2$-ovoids in polar spaces using a combination of spectral methods and the rank bound on $1$-ovoids. 
These bounds are weaker than the Ramsey bounds for partial $2$-ovoids, but they allow us to give a new simple proof of known non-existence results on $2$-ovoids \cite{BKLP2007}, and to obtain the following uniform consequence.

\begin{Theorem}\label{partial2ovoids_nonex}
Let $r \geq (1+\log_p(7))p+1$ for $p\geq 5$, or let $r\geq 6$ for $p=2,3$. 
Then a polar space of rank $r$ over a field of characteristic $p$ does not possess a $2$-ovoid.
\end{Theorem}

We then turn to constructions. 
In Section \ref{sec:constructions}, we give a general probabilistic construction of partial $m$-ovoids, which shows that for $m \geq cr$, with any constant $c > 1/2$, there is a partial $m$-ovoid of size $q^{\Omega(r)}$ in a polar space of rank $r$.
We write the result below with an explicit $c$ for convenience. 

\begin{Theorem} \label{probconstruction}
In a finite polar space of rank $r$, type $e$, over a finite field of order $q$, there exists a partial $m$-ovoid of size 
\[
q^{\Omega(r)}
\]
for every $m \geq 0.51r$.
\end{Theorem}

We then focus on a special case for explicit constructions.
Suppose we are in the symplectic space $\w(2r - 1, 2)$ defined over the vector space $\mathbb{F}_2^{2r}$. By an Oddtown argument, but also from the work of Dye \cite{Dye1992}, a partial ovoid of $\w(2r - 1, 2)$ has size at most $2r+1$.
For $m = 2$, the Ramsey number $R(3, 2r+2)$ yields the bound
$|\mathcal{O}_2| \le O(r^2/\log r)$. 
In Section \ref{connections}, 
we give an explicit construction of partial $2$-ovoids in $\w(2r - 1, 2)$ of size at least $(2r)^{1.12895}$, which greatly improves the trivial lower bound of $4r$ obtained by taking two disjoint partial $1$-ovoids of size $2r$. 
Note that \cite{Codenotti200} gives a better construction of size $8r - 4$, using a different terminology, and this was recently improved to roughly $(2r)^{1.088}$ in \cite{Golovnev2022}. 
Our construction improves on these bounds.

In $\w(2r - 1, 2)$, partial $m$-ovoids are also linked to a certain generalisation of the well-known Oddtown set families \cite{Berlekamp}. 
In Section~\ref{connections}, we state this generalised version of the Oddtown problem, and show its connection to the finite geometry problems of partial $m$-ovoids and $m$-nearly orthogonal sets. 
This allows us to translate constructions of triangle-free graphs with low complementary $2$-rank into partial $2$-ovoids and nearly orthogonal sets over the binary field.
One of our main contributions is a novel construction of such graphs using certain BCH codes.

\begin{Theorem}\label{thm:nearorth}
There exists an explicit $2$-nearly orthogonal set of size at least $t^{1.12895}$ in $\mathbb{F}_2^t$. 
Equivalently, there exists an explicit partial $2$-ovoid of $\w(t-1,2)$, $t$ even, of size at least $t^{1.12895}$.
\end{Theorem}

The construction uses triangle-free graphs associated to certain BCH codes whose complements have low $2$-rank. 
These graphs were used by Alon \cite{Alon95} to give explicit constructions for Ramsey graphs of triangles versus cliques, and more recently by Deng et al. \cite{DHWX2023} in the study of storage codes. 

In Section \ref{section:cap}, we show 
that for all $\varepsilon > 0$, there exists an $m>1$ such that for all large enough $r$,
$\w(2r-1, 2)$ possesses a partial $m$-ovoid of size at least
\[
(2r)^{\frac{\log(16)}{\log(10)}-\varepsilon}.
\]
The convergence is quite sharp and it leads to new constructions for various small values of $m$. 
Our construction is somewhat surprising from the viewpoint of finite geometry: a triangle-free
graph with low complementary $2$-rank is constructed, and then infinitely many 
graphs with low complementary $2$-rank are produced via the strong product, using
a result of \cite{BLS2024}. 
These graphs can be embedded in binary symplectic spaces
of a certain dimension because of their connection to \emph{nearly orthogonal sets} \cite{Chawin24}.

In this section, we also construct new triangle-free graphs using a binary projective cap that have low complementary rank over the reals. 
This improves the bounds in the recently introduced rank-Ramsey problem \cite{BLS2024} for triangles, which was motivated by the famous log-rank conjecture.
The best construction given in \cite{BLS2024} had rank $(3/8 + o(1))N$, whereas we prove the following.

\begin{Theorem}\label{cor:rank_ramsey_cap}
There exists a triangle-free graph $G$ on $N$ vertices with 
\[
\mathrm{rank}(A_G + I) = (1/4 + o(1))N
\]
over $\mathbb{R}$.
\end{Theorem}

We use the same construction to give improved lower bounds on partial $m$-ovoids in binary symplectic spaces for small values of $m$, and thus on $m$-nearly orthogonal sets over $\mathbb{F}_2^t$. 

\begin{Corollary}\label{cor:partovoidW}
For every fixed $m \geq 2$ and $t$ even and sufficiently large,
the binary symplectic space $\w(t-1, 2)$ has a partial $m$-ovoid of size at least 
\[t^{\beta_m - o(1)},\]
where 
\[\beta_m =  \log_{10}\left(\frac{32}{\left(32+144\cdot 2^{m+1}\right)^{1/(m+1)}}\right).\]
\end{Corollary}

This construction is better than the previous construction of partial $2$-ovoids for all $m \geq 28$.
In particular, for $t$ sufficiently large, we get that $\w(t,2)$ has a partial $28$-ovoid of size at least $t^{1.12936}$.

\section{Background on finite classical polar spaces}\label{background}

Throughout, we will use the symbol $q$ for a prime power $q:=p^h$, $p$ prime and $h$ is a positive integer, and we
will denote the finite field of order $q$ as $\F_q$. The vector space of dimension $n$ over
$\F_q$ will be written as $V(n,q)$, and $\pg(n-1,q)$ will denote the projective space with
underlying vector space $V(n,q)$. Let $f$ be a (reflexive) sesquilinear or quadratic form on
$V(n,q)$. The elements of the finite classical polar space $\p$ associated with $f$ are the
totally singular or totally isotropic subspaces of $\pg(n - 1,q)$ with relation to $f$, according to
whether $f$ is a quadratic or sesquilinear form. The Witt index of the form $f$ determines the vector space
dimension of the maximal subspaces contained in $\p$; the {\em rank} of $\p$ equals the
Witt index of its form. 
Hence, a finite classical polar space of rank $r$ embedded in $\pg(n - 1,q)$ has an underlying
form of Witt index $r$, and contains points, lines, \ldots, $(r-1)$-dimensional projective subspaces. The
subspaces of maximal dimension contained in $\p$ are called its \emph{generators}.
The (projective) dimension of generators will be one less than the Witt index.
Recall that a partial $m$-ovoid is a set of points that meets every generator in at most $m$ points, and it is called an $m$-ovoid if it meets every generator in exactly $m$ points. 
Also recall that for $m = 1$ these are simply called partial ovoids, and ovoids, respectively.

We will use projective notation for finite polar spaces so that they differ from the standard
notation for their collineation groups. For example, we will use the notation $\w(n-1,q)$ to denote
the symplectic polar space coming from the vector space $V(n,q)$, for $n$ even, equipped with a non-degenerate
alternating form. 
Here is a summary of the notation we will use for finite polar spaces, together
with their ranks and ovoid numbers, that is, the size of a putative ovoid in the polar space. 
In this table the parameter $e$ helps us compute various other combinatorial parameters of the polar space,
and essentially arises from knowing the number of points of the `rank 1' polar space of the same type
(which is $q^e+1$, where $q$ is the order of the defining field).
For example, a polar space of rank $r$ and type $e$ over a finite field $\mathbb{F}_q$ has $(q^r - 1)(q^{r + e - 1} + 1)/(q - 1)$ points and 
$\prod_{i = 1}^{r} (q^{i + e - 1} + 1)$ generators. The size of a putative ovoid, written as the Ovoid Number, is equal to $q^{r + e - 1} + 1$. 
A survey of ovoids in polar spaces can be found in \cite{DBKM}.

\begin{table}[H]
\caption{Notation for the finite classical polar spaces, together with their ovoid numbers.}\label{table:overview}
\begin{tabular}{lcccc}
\toprule
Polar Space&Notation&Ovoid Number& Type $e$ & Rank $r$\\
\midrule
Symplectic&$\w(n-1,q)$, $n$ even&$q^{n/2}+1$ & $1$ & $n/2$\\
Hermitian&$\h(n-1,q^2)$, $n$ odd&$q^{n}+1$ & $3/2$ & $(n - 1)/2$\\
Hermitian&$\h(n-1,q^2)$, $n$ even&$q^{n-1}+1$ & $1/2$ & $n/2$\\
Orthogonal, elliptic&$\mathsf{Q}^-(n-1,q)$, $n$ even&$q^{n/2}+1$ & $2$ & $(n - 2)/2$\\
Orthogonal, parabolic&$\mathsf{Q}(n-1,q)$, $n$ odd&$q^{(n-1)/2}+1$ & $1$ & $(n - 1)/2$\\
Orthogonal, hyperbolic&$\mathsf{Q}^+(n-1,q)$, $n$ even&$q^{n/2-1}+1$ & $0$ & $n/2$\\
\bottomrule
\end{tabular}
\end{table}

Let $\p$ be a finite polar space defined by a sesquilinear or quadratic form $f$, and let $X$ be a
point of the ambient projective space. Then $X^\perp$ is the set of projective points whose
coordinates are orthogonal to $X$ with respect to the form $f$. Note that when $f$ is a quadratic
form, it determines a (possibly degenerate when $q$ is even) bilinear form\footnote{When $f$ is a
quadratic form, $g(v,w) := f(v+w)-f(v)-f(w)$ is a bilinear form, symmetric when $q$ is odd, alternating when $q$ is even.} $g$, and two projective points
$X$ and $Y$ are orthogonal with relation to $f$ if, by definition, they are orthogonal with relation
to $g$.
A partial ovoid is thus a set of pairwise non-orthogonal points of a polar space. 
We refer to \cite{Hirschfeld1991} for further details on polar spaces and their combinatorial applications.

\subsection{Some \texorpdfstring{$p$}{p}-rank bounds for partial ovoids of finite polar spaces}\label{sec:pranks}

Our results in the next two sections rely on $p$-rank bounds for partial ovoids
which are summarised in this section. Here a \text{$p$-rank bound} is a bound that
corresponds to the $p$-rank of a matrix. In all the bounds below, either the point-hyperplane
incidence matrix of the underlying projective space (or a closely related matrix) give the bounds.

\begin{Theorem}[Folklore, cf.~\cite{AS2011,BM1995}] \label{thm:symp_bnd}
  Let $\mathcal{O}$ be a partial ovoid in a symplectic space naturally embedded
  in a vector space of even
  dimension $n \geq 4$ over a field with $q = p^h$, where $p$ is prime, elements.
  Then 
  \begin{align*}
    |\mathcal{O}| \leq \binom{p+n-1}{p-1}^h +1.
  \end{align*}
\end{Theorem}

\begin{Theorem}[{\cite[Theorem 1.6]{AS2011}}]\label{thm:quadric_bnd}
  Let $\mathcal{O}$ be a partial ovoid in a quadric in a vector space of 
  dimension $n \geq 4$ over a field with $q = p^h$, where $p$ is prime, elements.
  Then 
  \begin{enumerate}[(a)]
   \item If $p=2$ and $n$ even, then $|\mathcal{O}| \leq n^h+1$.
   \item If $p=2$ and $n$ odd, then $|\mathcal{O}| \leq (n-1)^h+1$.
   \item If $p>2$, then 
   \[
      |\mathcal{O}| \leq \left[ \binom{n+p-2}{p-1} - \binom{n+p-4}{p-3} \right]^h + 1.
   \]
   \item If $p>2$ and there exists a positive integer $u$ satisfying
   \[
     u+1 \equiv n \pmod{2} \text{ and } n-3 \leq up \leq n+p-5,
   \]
   then 
   \[
    |\mathcal{O}| \leq \left[ \binom{n+p-2}{p-1} - \binom{n+p-4}{p-3} - \binom{up+2}{n-1} + \binom{up}{n-1} \right]^h + 1.
   \]
  \end{enumerate}
\end{Theorem}

\begin{Theorem}[{\cite[Theorem 1.8]{AS2011}}]\label{thm:herm_bnd}
  Let $\mathcal{O}$ be a partial ovoid in a Hermitian polar space in a vector space of 
  dimension $n \geq 4$ over a field with $q = p^{2h}$, where $p$ is prime, elements.
  \begin{enumerate}[(a)]
   \item  We have 
   \[
      |\mathcal{O}| \leq \left[ \binom{n+p-2}{p-1}^2 - \binom{n+p-3}{p-2}^2 \right]^h + 1.
   \]
   \item If there exists a positive integer $u$ satisfying 
   \[
     n-2 \leq up \leq n+p-4,
   \]
   then 
   \[
    |\mathcal{O}| \leq \left[ \binom{n+p-2}{p-1}^2 - \binom{n+p-3}{p-2}^2 - \binom{up+1}{n-1}^2 + \binom{up}{n-1}^2 \right]^h + 1.
   \]
  \end{enumerate}
\end{Theorem}

\section{A Ramsey bound on partial \texorpdfstring{$m$}{m}-ovoids}
\label{sec:ramsey_part}

A classic result of Erd\H{o}s and Szekeres~\cite{es1935} states that $R(s, t) \leq \binom{s+t-2}{s-1}$, 
where $R(s, t)$ denotes the smallest number of vertices $v$ that
for any graph $\Gamma$ on $v$ vertices, either $\Gamma$ contains 
a clique of size $s$ or a coclique of size $t$. 
For fixed $s$, this upper bound is roughly $O(t^{s - 1})$, which was improved to 
\[R(s, t) \leq c_s \frac{t^{s - 1}}{(\log t)^{s - 2}}\]
by  Ajtai, Koml\'os and Szemer\'edi \cite{AKS1980}. 
For $s = t$, a recent breakthrough of Campos et al.~\cite{campos2023} showed that 
\[R(t,t) \leq (4 - \epsilon)^t,\]
thus making an exponential improvement over the classical upper bound of $\binom{2t - 2}{t - 1} \approx (1 + o(1)) \frac{4^t}{\sqrt{\pi t}}$.

We show that partial $m$-ovoids can give lower bounds for Ramsey numbers, or equivalently, upper bounds on Ramsey numbers imply upper bounds on the size of a partial $m$-ovoid.

\begin{proof}[Proof of Theorem \ref{Ramseypartialmovoid}]
Let $\Gamma$ be the collinearity graph of $\mathcal{P}$. 
Then $\Gamma$ does not contain a coclique of size $k + 1$, since $k$ is the maximum size of a partial ovoid. 
The induced subgraph $H$ on $\mathcal{O}_m$ does not contain a clique of size $m + 1$, as any maximal clique of $\Gamma$ corresponds to a generator of the polar space and hence meets $\mathcal{O}_m$ in at most $m$ points. 
Therefore $H$ is a graph on $|\mathcal{O}_m|$ vertices with no coclique of size $k + 1$ and no clique of size $m+1$, thus proving the bound. 
\end{proof}

\begin{proof}[Proof of Corollary \ref{Ramseypartialmovoid_corollary}]
    The size of an $m$-ovoid in a polar space of rank $r$ and type $e$ over $\mathbb{F}_q$, with $q = p^h$, is $m(q^{r + e - 1} + 1) = m(p^{hr + he - h} + 1)$.
    As $p \geq 2$ and $m \geq 2$, we can assume that $r \geq 7$.
    We want to apply Theorem \ref{Ramseypartialmovoid} with $R(m+1, k+1) \leq \binom{m+k}{m} \leq k^m$,
    where $k = \binom{p+n-1}{p-1}^h {{+ 1}}$. Here we used that $n \geq 2r \geq 14$ and $m \geq 2$.
    Hence, it suffices to show that if $r \geq m(p-1) \log_p(2r+2+p)+1$, then
    \begin{align*}
      m \cdot p^{hr + he - h} \geq k^m.
    \end{align*}
      Since
    \[
      \binom{p+n-1}{p-1} + 1\leq (n+p)^{p-1}\leq (2r+2+p)^{p-1},
    \]
    it suffices to show that
    \[
      m \cdot p^{hr+he-h} \geq (2r+2+p)^{mh(p-1)}.
    \]
    Taking logarithms, this follows from
    \[
      \log_p(m)+hr+he-h \geq mh(p-1)\log_p(2r+2+p).
    \]
    Note that 
    \[
      \log_p(m)+hr+he-h=\log_p(m)+h(r+e-1)\geq h(r-1) 
    \]
    so it is enough to prove
    \[
      r-1 \geq m(p-1) \log_p(2r+2+p).  \qedhere 
    \]
\end{proof}

We note that the results of \cite{BKLP2007}, \cite{DeBMS}, and \cite{HamiltonMathon}  on the non-existence of $m$-ovoids only apply to
the polar spaces $\q^-(2r+1,q)$, $\w(2r-1,q)$, and $\h(2r,q^2)$. In contrast, the bound in Corollary \ref{Ramseypartialmovoid_corollary} applies to all finite classical polar spaces, thus proving the non-existence of $m$-ovoids in all polar spaces of large rank.

\section{Spectral bounds on partial \texorpdfstring{$2$}{2}-ovoids and non-existence of \texorpdfstring{$2$}{2}-ovoids}\label{sec:spectralbounds}
In this section we show the bounds on partial $2$-ovoids that one can obtain by using spectral arguments and compare them to results based on Ramsey numbers. 
We observe that while the bounds on partial $2$-ovoids are worse, they allow us to show the non-existence of $2$-ovoids in a wider class of polar spaces. 

Let $P$ be a point of a finite polar space $\mathcal{P}$. We can project the elements of $\mathcal{P}$ incident with $P$
to (abstractly) obtain a finite polar space with the same type as $\mathcal{P}$, but with rank one less than
that of $\mathcal{P}$. For instance, each line of $\mathcal{P}$ incident with $P$ maps to a point
of the \emph{quotient} polar space. By using the non-existence of ovoids in the quotient polar space,
we can rule out the existence of 2-ovoids in the original polar space.
Recall the following non-existence result. 

\begin{Theorem}[{\cite[Theorem 14]{BKLP2007}}]\label{BKLP_Theorem14}\samepage
The following polar spaces do not admit $2$-ovoids:
\begin{enumerate}[(i)]
\item $\w(2r-1,q)$ for $r>2$ and $q$ odd; 
\item $\q^-(2r+1,q)$ for $r>2$;
\item $\h(2r,q^2)$ for $r>2$;
\item $\q(2r,q)$ for $r>4$.
\end{enumerate}
\end{Theorem}

We say that a graph, neither complete nor edge-less, is \emph{strongly regular}
with parameters $(v, k, \lambda, \mu)$ if it has $v$ vertices, degree $k$,
each pair of adjacent vertices $x,y$ has precisely $\lambda$ common neighbours,
and each pair of distinct nonadjacent vertices $x,y$ have precisely $\mu$
common neighbours. The eigenvalues of a strongly regular graph 
are $k \geq e^+ \geq 0 > e^-$ (cf.~ \cite{BvM}).
The collinearity graphs on the points of a polar space 
are strongly regular graphs. We will use the following well-known result (see \cite[Prop.\ 1.1.6]{BvM}).

\begin{Lemma}\label{lem:ind_gr}
 The induced subgraph on a vertex set $Y$ of a $k$-regular 
 graph $\Gamma$ of order $v$ with smallest eigenvalue $s$ 
 has average degree at least $|Y| \frac{k-s}{v} + s$.
\end{Lemma}

Now we can use both the information of the strongly regular graph, and the size of partial ovoids
of a quotient polar space, to bound the size of a partial 2-ovoid.

\begin{Theorem}\label{partial2ovoids} \label{W_2_ovoid}
  Let $r\ge 3$ and let $\mathcal{O}_2$ be a partial $2$-ovoid of a finite polar space $\mathcal{P}$ of rank $r$
  and with defining field of order $q$.
  Then $|\mathcal{O}_2|\le q \cdot b(n-2, q)+ \onum_r$ where $\onum_r$ is the ovoid number of $\mathcal{P}$,
  and $b(n', q)$ is the maximal size of a partial ovoid in a polar space 
  of the same type as $\cP$ in an $n'$-dimensional vector space over the field with $q$ elements.
\end{Theorem}

\begin{proof}
First, the collinearity graph $\Gamma$ on the points of $\mathcal{P}$
is a strongly regular graph. Let $\onum_r$ be the ovoid number of $\mathcal{P}$.
Then the parameters of $\Gamma$ (see \cite[\S2.5.2]{BvM})
are
\begin{align*}
  & v = \onum_r \frac{q^r-1}{q-1}, && e^+ = q^{r-1}-1,\\
  & k = q\frac{q^{r-1}-1}{q-1}\onum_{r-1}, && e^- = -\onum_{r-1}.
\end{align*}
By Lemma \ref{lem:ind_gr}, we find an element $P \in \mathcal{O}_2$ adjacent 
to at least $M := |\mathcal{O}_2| \frac{k-e^-}{v} + e^-$ vertices.
Hence, as no line contains more than $2$ points, 
$\mathcal{O} := \{ \< P, Q \>/P: Q \in {\mathcal{O}_2 \cap P^\perp, Q \neq P}\}$ has size at least $M$.
Since $\mathcal{O}$ is a partial $1$-ovoid of $P^\perp/P$,
we have 
  \[
    |\mathcal{O}_2| \frac{k-e^-}{v} \leq b(n-2, q) - e^-.
  \]

Thus,
\begin{align*}
|\mathcal{O}_2| &\le \frac{v}{k-e^-}\left(b(n-2, q) - e^-\right)\\
&=\frac{\onum_r \frac{q^r-1}{q-1}}{ (q\frac{q^{r-1}-1}{q-1}+1)\onum_{r-1}}
\left(b(n-2, q) +\onum_{r-1}\right)\\
&=\frac{\onum_r}{ \onum_{r-1}}b(n-2, q)+ \onum_r.
\end{align*}
Now $\onum_r = q\onum_{r-1}-q+1$ and so $|\mathcal{O}_2| \le q \cdot b(n-2, q)+ \onum_r$.
\end{proof}

Together with Theorem \ref{thm:symp_bnd}, Theorem \ref{thm:quadric_bnd},
and Theorem \ref{thm:herm_bnd}, we have the following:

\begin{Corollary}\label{partial2ovoids_corollary}
  Let $r\ge 3$ and let $\mathcal{O}_2$ be a partial $2$-ovoid of a finite polar space $\mathcal{P}$ of rank $r$
  and with defining field of order $q$. Then
  \[
      |\mathcal{O}_2| \le
  \begin{cases}
 q^{r} + q+1+ q \left[ \binom{p+2r-3}{p-1} \right]^h  & \text{if }\mathcal{P}\cong\w(2r-1,q),\\
 q^{r-\eps}+q+1 + q \left[ \binom{n+p-4}{p-1} - \binom{n+p-6}{p-3} \right]^h 
 & \text{if }\mathcal{P}\cong\q^\eps(2r-1,q),\\
 q^{2r+ 1}+q^2+1+ q^2 \left[ \binom{n+p-4}{p-1}^2 - \binom{n+p-5}{p-2}^2 \right]^h 
  & \text{if }\mathcal{P}\cong\h(n-1,q^2), n\text{ odd},\\
   q^{2r- 1}+q^2+1 + q^2 \left[ \binom{n+p-4}{p-1}^2 - \binom{n+p-5}{p-2}^2 \right]^h 
  & \text{if }\mathcal{P}\cong\h(n-1,q^2), n\text{ even}.
  \end{cases}
    \]
\end{Corollary}

Note that Theorem \ref{BKLP_Theorem14}
is also implied by Corollary \ref{partial2ovoids_corollary}.

\begin{proof}[Proof of Theorem \ref{partial2ovoids_nonex}]
  The size of a $2$-ovoid in a polar space of rank $r$ and type $e$ over $\Ff_q$,
  with $q = p^h$, is $2(q^{r+e-1}+1) = (q^{r+e-1} + 1) + (p^{h(r+e-1)} + 1)$; see Table \ref{table:overview}.
  Note that the first summand also occurs in the bounds in Corollary \ref{partial2ovoids_corollary},
  so we only need to show that $p^{h(r+e-1)} + 1$ is larger than $q \binom{p+2r-3}{p-1}^h + q-1$.

  First we show the claim for $p \geq 3$.
  Using $p-1 \geq 2$, we have
  \begin{align*}
    q \binom{p+2r-3}{p-1}^h + q - 1 \leq p^h \cdot (2r+e+p)^{h(p-1)}.
  \end{align*}
  We want to show that
  \begin{align*}
    p^{r+e-1} > p (2r+e+p)^{p-1}
  \end{align*}
  as then $2$-ovoids cannot exist.
  This is implied by 
  \begin{align*}
      r+e-2 > (p-1)\log_p(2r+e+p).
  \end{align*}
  Hence, it suffices to show that 
  \begin{align*}
      r-2 > (p-1)\log_p(2r+p+2).
  \end{align*}
  For $r = (1+\log_p(7)) p +1$, we find
  \begin{align*}
    \log_p(2r+p+2) = \log_p((3+2\log_p(7))p + 4) \leq \log_p((4+2\log_p(7))p) = 1 + \log_p(4 + 2\log_p(7)).
  \end{align*}
  As $2\log_p(7) < 3$ for $p \geq 5$, this shows $r-2 > (p-1) \log_p(2r+p+2)$.
  For $p=2,3$, it can be verified that $r \geq 6$ suffices.
\end{proof}

Corollary \ref{Ramseypartialmovoid_corollary} gives the
condition of the form $r \geq (\log_p (2r+2+p)^2)(p-1)+1$
which is worse than Theorem \ref{partial2ovoids_nonex} for $p \geq 3$.

 \section{Constructions} 
 
In this section we give a general probabilistic construction of partial $m$-ovoids, when $m$ grows at least linearly with respect to the rank $r$ of the polar space, and an explicit construction for the special case of $m = 2$ in the binary symplectic polar space, via an equivalence with \emph{nearly orthogonal sets}. 
We also give constructions of partial $m$-ovoids for $m > 2$ that have a probabilistic component in them. 
\subsection{Probabilistic construction}\label{sec:constructions}

In this section we give a general probabilistic construction of partial $m$-ovoids, when $m$ grows at least linearly with respect to the rank $r$ of the polar space.

\begin{proof}[Proof of Theorem \ref{probconstruction}]
    Let $\mathcal{P}$ be a finite polar space of rank $r$ and type $e$ over the finite field $\mathbb{F}_q$. 
Let $N = \theta_r (q^{r + e - 1} + 1)$ be the total number of points in $\mathcal{P}$, where $\theta_r = (q^r - 1)/(q - 1) \leq q^r$ and let $M = \prod_{i = 1}^r (q^{i + e - 1} + 1)$ be the total number of generators. 
Pick a random subset $S$ of $\mathcal{P}$ by picking each point independently with probability $\rho$, where the exact value of $\rho$ will be specified later.
Thus, $S$ is the outcome of a Bernoulli process, that is a sequence of $N$ independent identically distributed Bernoulli trials with parameter $\rho$.

For a fixed generator, the probability that $S$ contains at least $m + 1$ points from it is at most
\[\binom{\theta_r}{m + 1}\rho^{m+1}.\]
Therefore, by the union bound (for probabilities) over all generators, if we have 
\[M \binom{\theta_r}{m+1} \rho^{m+1} < 1/2,\] then with probability $> 1/2$, $S$ does not contain $m+1$ points from any generator, that is, it is a partial $m$-ovoid.
Let $A = r(r + 1)/2 + r(e-1)$.
Since $M = \prod_{i = 1}^r (q^{i + e - 1} + 1) \leq \prod_{i = 1}^r 2q^{i + e - 1} = 2^r q^{r(r + 1)/2 + r(e-1)} \leq q^{A + r}$ and $\binom{\theta_r}{m+1} < q^{r(m+1)}$, 
taking \[\rho = \frac12 q^{-r - \frac{A + r}{m+1}},\] we get $M \binom{\theta_r}{m+1} \rho^{m+1} \leq q^{A+r}q^{r(m+1)} \left( \frac12 q^{-r - \frac{A + r}{m+1}} \right)^{m+1} = 2^{-(m+1)} < 1/2$.
Therefore, with this value of $\rho$ $\mathrm{Pr}(S \text{ is a partial }m\text{-ovoid}) > 1/2$.

Let $X=|S|$. Since $N\geq q^{2r+e-2}$, we have \[ \mathbb{E}[X]=\rho N \geq \frac12 q^{r + e - 2 - \frac{A + r}{m+1}}. \] 
For $m\geq 0.51r$, the exponent is at least \[ \left(1-\frac{1}{1.02}\right)r-O(1)=\Omega(r), \] and hence $\mathbb{E}[X]=q^{\Omega(r)}$. By a standard Chernoff bound, $\Pr(X<\mathbb{E}[X]/2)=o(1)$. 
Therefore, with positive probability, $S$ is a partial $m$-ovoid and has size at least $\mathbb{E}[X]/2=q^{\Omega(r)}$. This proves the result.
\end{proof}

\subsection{Oddtown, nearly orthogonal sets, and triangle-free graphs}\label{connections}

A family $\mathcal{F}$ of subsets of a $t$-element set is an \emph{Oddtown} if $|S|$ is odd for every 
$S\in \mathcal{F}$, and $|S\cap T|$ is even for every pair of distinct elements $S,T\in\mathcal{F}$.
It was shown by Berlekamp \cite{Berlekamp}, that the size of an Oddtown is at most $t$ using a simple linear algebraic argument. There is a direct connection between partial $m$-ovoids of binary symplectic spaces and 
similar families of subsets. Let $V$ be the $2r$-dimensional vector subspace of $\mathbb{F}_2^{2r + 1}$ consisting of all even weight vectors, i.e., $V = \{(x_1, \dots, x_{2r + 1}) : \sum_{i = 1}^{2r + 1} x_i = 0\}$. Then $V$ along with the ``standard'' dot product, $(x_1, \dots, x_{2r + 1}) \cdot (y_1, \dots, y_{2r+1}) = \sum_{i = 1}^{2r+1} x_i y_i$ gives rise to a symplectic polar space isomorphic to $\w(2r - 1, 2)$. 
Let $\mathcal{O}$ be a partial $m$-ovoid in this model of $\w(2r - 1, 2)$. 
Then the vectors representing elements of $\mathcal{O}$ have the following properties:
\begin{enumerate}[(i)]
    \item For all $x \in \mathcal{O}$, $x \cdot x = 0$. 
    \item For any $(m + 1)$ distinct points $x_1, \dots, x_{m + 1} \in \mathcal{O}$ there exists $i \neq j$ such that $x_i \cdot x_j = 1$. 
\end{enumerate}
If we think of vectors here as characteristic vectors of subsets of $\{1,\ldots, 2r+1\}$, then in terms of set families, we have a family $\mathcal{F} \subseteq 2^{[2r + 1]}$ such that
for all $S \in \mathcal{F}$, $|S|$ is even, and 
for all distinct $S_1, \dots, S_{m +1} \in \mathcal{F}$, there exists an $i \neq j$ such that $|S_i \cap S_j|$ is odd.
By taking complements of these sets, we obtain the following generalisation of an Oddtown:

\begin{Definition}[Generalised Oddtown]
A family $\mathcal{G}$ of subsets of a $t$-element set is a \emph{generalised Oddtown}
with parameter $m$, if:
\begin{enumerate}[(i)]
    \item for all $S \in \mathcal{G}$, $|S|$ is odd;
    \item for all distinct $S_1, \dots, S_{m +1} \in \mathcal{G}$, there exists an $i \neq j$ such that $|S_i \cap S_j|$ is even.
\end{enumerate}
When $m=1$, we have an Oddtown.
\end{Definition}

Let $\Gamma_t$ be the graph whose vertices are the proper subsets of $\{1,\ldots, t\}$ of odd cardinality,
and two such subsets are adjacent if they are not equal and have intersection of odd cardinality.
For $t$ odd, $\Gamma_t$ is isomorphic to the collinearity graph of the symplectic space $\w(t-2,2)$,
and for $t$ even, $\Gamma_t$ is isomorphic to the induced subgraph of the collinearity graph, on the points of
$\w(t-1,2)$ not lying in a fixed hyperplane. Therefore, the clique number of $\Gamma_t$ is $2^{(t-1)/2} - 1$, for $t$ odd.
A generalised Oddtown with parameter $m$ translates to a subset of the vertices of $\Gamma_t$
not containing an $(m+1)$-clique. So if $m>2^{(t-1)/2}$, then the set of all subsets of $\{1,\ldots,t\}$ of odd
cardinality is a generalised Oddtown.
 
To avoid unnecessary case distinctions depending on the parity of $t$, we will
assume that $t$ is odd and write $\w(t-2, 2)$. If $t$ is even, then everything we
write holds for the graph obtained above from $\w(t-1, 2)$.
Recall that a subset $\cF$ of vectors of $\F^t$ is called \emph{$m$-nearly orthogonal} if its 
elements are not self-orthogonal
and every subset of $m+1$ members includes an orthogonal pair.
For $\F = \F_2$, it is clear from the
discussion above that this is equivalent to our generalised Oddtown problem as well as to the maximal
size of a partial $m$-ovoid in $\w(t-2, 2)$. 
Thus, the results of \cite{Golovnev2022} on nearly orthogonal sets imply
 that there exists a partial $2$-ovoid of $\w(t-2, 2)$
of size $t^{1+\delta}$ for some constant $\delta > 0$.  More recently, it has been shown in Theorem 1.1 of \cite{Chawin24}
that there exists a constant $\delta > 0$ such that for all integers $t \geq m \geq 2$,
there exists a partial $m$-ovoid of $\w(t-2, 2)$ if $t$ is odd, and of
$\underline{\w(t-1,2)}$ if $t$ is even, of size at least
\[
 t^{\delta m/\log m}.
\]
Note that this does not imply anything for small $m$.
In all of the aforementioned results, $\delta$ can be calculated without too much effort
from the respective proofs. For instance, the proof of \cite[Theorem 1.1]{Chawin24} shows that we find a partial
$800$-ovoid of size $100^{1.25\ell-C}$ in $\w(100^\ell-1, 2)$.

Here we give an explicit construction for $2$-nearly orthogonal sets, and hence partial $2$-ovoids that improves these previous results. 
We follow the proof of Theorem 1.11 in \cite{Golovnev2022}, but starting from a triangle-free graph related to BCH codes instead of the generalized Kneser graphs. 
This graph was used by Alon in 1995 \cite{Alon95} to give an explicit construction for Ramsey graphs of triangles vs cliques, and more recently by Deng et al. \cite{DHWX2023} to solve a problem related to storage codes. 
We first state the following lemma of Lempel \cite{lempel} that's crucial for our result. 

\begin{Lemma}
\label{lem:Lempel}
Let $M$ be an $n \times n$ symmetric matrix over the binary field $\mathbb{F}_2$ with at least one nonzero diagonal entry and rank $r$. Then, there exists an $n \times r$ matrix $B$ over $\mathbb{F}_2$ satisfying $M = B B^T$.
\end{Lemma}

\begin{proof}[Proof of Theorem \ref{thm:nearorth}]
  Consider the Cayley graph $\Gamma = \mathrm{Cay}(G, S)$ with $G = (\mathbb{F}_{2^h}^2, +)$ and $S = \{(a, a^3) : a \in \Ff_{2^h} \setminus \{ 0 \} \}$.
  This graph is triangle-free because there are no three $a, b, c \in \Ff_{2^h} \setminus \{0\}$ with $a + b + c = a^3 + b^3 + c^3 = 0$. 

    Let $M = J - A' = A + I$, where $J$ is the all-ones matrix, $A$ is the adjacency matrix of $\Gamma$, and $A'$ is the adjacency matrix of the complement $\overline{\Gamma}$. 
  In Theorem 2 of \cite{DHWX2023} it is shown that $M$ has $2$-rank at most
  \[
    t := \left\lfloor \frac{1+\sqrt{2}}{2} ( 2 + \sqrt{2} )^h \right\rfloor.
  \]
  By Lemma~\ref{lem:Lempel}, $M = B B^T$ for some $n \times t$ matrix $B$, where $n = 2^{2h}$. 
  The row vectors of $B$ give us $n$ vectors in $\mathbb{F}_2^t$ that are not self-orthogonal.
   As $\Gamma$ is triangle-free, for any three of these vectors, at least two of them must be orthogonal to each other with respect to the standard dot product. 
  Hence, we obtain a $2$-nearly orthogonal set of the asserted size (note: $\log(4)/\log(2+\sqrt{2})\approx 1.12895$).
\end{proof}

\subsection{Triangle-free graphs from projective caps}
\label{section:cap}
A \emph{cap} in $\mathrm{PG}(n - 1, 2)$ is a set $S$ of points such that no three of them are collinear. 
In terms of the underlying vector space, this corresponds to a set of non-zero vectors such that no three of them sum to the zero vector. 
Given such a set, we can construct the Cayley graph of the abelian group $\mathbb{F}_2^n$ with $S$ as its generating set. 
This Cayley graph is clearly triangle-free.
From the character theory of finite abelian groups (see for example \cite[Chapter 11]{GodsilMeagher}), it also follows that $\lambda$ is an eigenvalue of this graph if and only if there is a vector $v \in \mathbb{F}_2^n$ for which $\sum_{u \in S} (-1)^{u \cdot v} = \lambda$. 
The zero vector corresponds to the eigenvalue $|S|$, which is the degree of the graph, and every other vector $v$ gives rise to a hyperplane $H_v$ defined by $\{x \mid v \cdot x = 0\}$ that meets $S$ in exactly $(|S| + \lambda)/2$ points. 
In conclusion, a cap $S$ of $\mathrm{PG}(n - 1, 2)$ gives rise to a $|S|$-regular triangle-free graph which has an eigenvalue $\lambda \neq |S|$ with multiplicity $m$ if and only if there are $m$ hyperplanes that intersect $S$ in exactly $(|S| + \lambda)/2$ points. 

For the rest of this section we let $S$ be the following cap in $\mathrm{PG}(n - 1, 2)$. 
Let $M$ be a codimension-$1$ subspace, and let $W$ be a codimension-$2$ subspace of $\mathrm{PG}(n - 1, 2)$ contained in $M$. Define $Z \coloneqq M \setminus W$. Let $P$ be a point in $Z$. Consider a line $\ell$ in $\mathrm{PG}(n-1, 2)$ through $P$ such that $\ell \cap M = \{P\}$.
We define 
\[
S \coloneqq (Z \cup \ell) \setminus \{P\},
\]
 and note that $|S|=|Z| + 1 = 2^{n-2}+1$ (see Figure~\ref{fig:cap}). 

\begin{figure}[h!]
\begin{center}
\begin{tikzpicture}[x={(1cm,0.4cm)}, y={(8mm, -3mm)}, z={(0cm,1cm)}, line cap=round, line join=round]
    \coordinate (x1) at (-2,2,3);
    \coordinate (x2) at (2,2,5);
    \coordinate (x3) at (2,-2,5);
    \coordinate (x4) at (-2,-2,3);

    \coordinate (P) at (-2.5, -0.25 , 4.1);

    \coordinate (Q) at ($(P) + (0,0,2)$);  
    \coordinate (R) at ($(P) + (0,0,4)$);  

    \path[draw=black, fill=black!20, thick, opacity=0.8] (x1) -- (x2) -- (x3) -- (x4) -- cycle;
    \node at ($(x3)!0.5!(x4) + (1, 0, 1.2)$) {$M$};

    \coordinate (y1) at (-1, 1, 3.5);
    \coordinate (y2) at (1, 1, 4.5);
    \coordinate (y3) at (1, -1, 4.5);
    \coordinate (y4) at (-1, -1, 3.5);
    
    \path[draw=black, fill=red!30, thick, opacity=0.8] (y1) -- (y2) -- (y3) -- (y4) -- cycle;
    \node at ($(y3)!0.5!(y4) + (0.2, 0, -0.5)$) {$W$};

    \draw[thick] (P) -- (R);

    \draw[fill=black] (P) circle (2pt) node[below right] {$P$};

    \node at ($(x3)!0.5!(x4) + (-0.7, 0, 1.7)$) {$\ell$};

\end{tikzpicture}
\end{center}
\caption{The cap $S$ in $\mathrm{PG}(n - 1, 2)$ is the complement of $W$ in $M$,
together with the points of $\ell$, minus the point $P$.}
\label{fig:cap}
\end{figure}

The construction gives us a cap as any line in $\mathrm{PG}(n - 1, 2)$ consists of three points and at most two of them can be contained in $S$. 

\begin{Remark}
Choosing $M = \{(x_0, \dots, x_{n - 1}) \mid x_0 = 0 \}$, $W = \{(x_0, \dots, x_{n-1}) \mid x_0=0, x_1=0 \}$, and $\ell$ to be the line spanned by $(0, 1, 0, \dots, 0)$ and $(1, 1, \dots, 1)$ gives us the construction in \cite{davy2, gabi}.
\end{Remark}

We now compute the possible intersection sizes of the hyperplanes with $S$, and the corresponding multiplicities. 

\begin{Lemma}
    Every hyperplane in $\mathrm{PG}(n - 1, 2)$ intersects $S$ in $1$, $\frac{|S|-3}{2}$, $\frac{|S|+1}{2}$, or $|S|-2$ points, with the number of hyperplanes intersecting in these many points is equal to $2$, $2^{n - 2} - 1$, $3 (2^{n - 2} - 1)$, and $1$, respectively.
\end{Lemma}

\begin{proof}
Consider a hyperplane $H$ in $\pg(n-1, 2)$. Since $S$ consists of points in $M$ and $\ell$, the intersection of $H$ with $S$ depends on how $H$ intersects $M$ and $\ell$. 
Recall that $M$ is a hyperplane in $\pg(n-1, 2)$,  $W$ is a codimension-2 subspace in $\pg(n - 1, 2)$ contained in $M$, and $Z=M \setminus W$.
We will repeatedly use the following argument: if the hyperplane $H$ is not equal to $M$, then it meets $M$ in a codimension-2 subspace, which is either equal to $W$ or it meets $Z$ in exactly $|W| - |W \cap H| = (2^{n - 2} - 1) - (2^{n - 3} - 1) = 2^{n - 3}$ points.
We proceed with three cases depending on whether $\ell \subseteq H$ or $|H \cap \ell| = 1$. 

\subsection*{Case 1: $\ell \subseteq H$}
The number of hyperplanes containing $\ell$ is $2^{n-2}-1$.
Since $P \in H \cap M$ and $P \notin W$, we have $H \cap M \neq W$. 
Therefore, $H$ meets $Z$ in $2^{n-3}$ points and that  includes the point $P$. Thus, $|H \cap S| = (2^{n - 3} - 1) + 2 = \frac{|S|+1}{2}$.

\subsection*{Case 2: $P \in H$ and $\ell \not\subseteq H$}
The number of hyperplanes passing through $P$ is $2^{n-1}-1$, and the number of hyperplanes passing through the line $\ell$ is $2^{n-2}-1$. Therefore, the number of hyperplanes passing through $P$, but not through $\ell$, is  $2^{n-2}$. 
One such hyperplane is $H = M$ and then $|H \cap S| = |Z| - 1 = |S| - 2$. Otherwise, $H \cap M \neq W$ because $P \in H$, and there are $2^{n-2}-1$ such hyperplanes. In these cases, $H$ meets $Z$ in $2^{n-3}$ points, including $P$, so $|H \cap S| = 2^{n - 3} - 1 = \frac{|S|-3}{2}$.

\subsection*{Case 3: $P \notin H$}
Let $\ell \cap H = \{Q\}$. 
Just like the previous case, the number of hyperplanes passing through $Q$, but not through $\ell$, is equal to $2^{n - 2}$. 
If $H \cap M = W$, then $H$ must be the span of $W$ and $Q$, and thus $|H \cap S| = 1$. 
Otherwise, $H \cap M \neq W$, and there are $2^{n-2}-1$ such hyperplanes. 
In these cases, $H$ meets $Z$ in $2^{n-3}$ points, so $|H \cap S| = 2^{n - 3}  + 1 = \frac{|S|+1}{2}$.
Note that there are two points in $\ell \setminus \{P\}$, and thus in this case, we get two hyperplanes that meet $S$ in $1$ point, and $2(2^{n- 2} - 1)$ hyperplanes that meet $S$ in $\frac{|S| + 1}{2}$ points. 
\end{proof}

\begin{Corollary}
    The spectrum of the Cayley graph of $\mathbb{F}_2^n$ with generating set equal to $S$ is 
    \[(2^{n - 2} + 1)^{(1)}, (2^{n - 2} - 3)^{(1)}, 1^{3(2^{n - 2} - 1)}, (-3)^{(2^{n - 2} - 1)}, (1 - 2^{n - 2})^{(2)}\]
\end{Corollary}
\begin{proof}
    The eigenvalue $2^{n - 2} + 1 = |S|$ corresponds to the degree of the graph. 
    The remaining eigenvalues are equal to $2|H \cap S| - |S|$, as $H$ goes through all hyperplanes of $\pg(n - 1, 2)$. 
\end{proof}

\begin{proof}[Proof of Theorem \ref{cor:rank_ramsey_cap}]
Following \cite[Corollary 4.5]{BLS2024}, we take the Kronecker product (also known as the tensor product) of this Cayley graph with the complete graph $K_\ell$, to get a triangle-free graph $G$ on $N = 2^n \ell$ vertices such that the multiplicity of the eigenvalue $-1$ is equal to $(\ell - 1) 3 (2^{n - 2} - 1)$. 
This gives us that the real rank of $A_G + I$ is equal to $(2^n \ell + 2^n \cdot 3 + 12 \ell - 12)/4$, which is $(1/4 + o(1)) N$ when $\ell, n \rightarrow \infty$. 
\end{proof}

We now use this graph to give a construction of large partial $m$-ovoids in $\w(2r - 1, 2)$ for $m > 2$. 
In \cite[Theorem 4.6]{BLS2024}  we find the following result, where we note that their proof works for any field,
not only the reals:

\begin{Theorem}\label{thm:BLS}
  Let $m > 1$ be a positive integer.
  Let $\Gamma$ be a graph on $n$ vertices and $r := \mathrm{rank}_\Ff(A_\Gamma+I)$
  for some field $\Ff$. Then, for every large enough $h$, there exists a graph $\Delta$
  such that:
  \begin{enumerate}
   \item the clique number of $\Delta$ is at most $m$,
   \item $\mathrm{rank}_\Ff(A_\Delta+I) \leq r^h$,
   \item the graph $\Delta$ has
   \[
    \frac{1}{m} \left( \frac{n}{\left( \sum_{t=1}^m t^{m+1} a_t \right)^{\frac{1}{m+1}} } \right)^h - 1
   \]
  vertices. Here $a_i$ denotes the number of cliques of size $i$ in $\Gamma$.
  \end{enumerate}
\end{Theorem}

Using this result and our triangle-free graph with low complementary rank, we prove the following. 

\begin{Proposition}\label{low2rank}
For every fixed integer $m\geq 2$ and all sufficiently large $t$, there exists a
$K_{m+1}$-free graph with complementary $2$-rank exactly $t$ and at least
\[
t^{\beta_m-o(1)}
\]
vertices, where
\[
\beta_m :=
\log_{10}\left(
\frac{32}{\left(32+144\cdot 2^{m+1}\right)^{1/(m+1)}}
\right).
\]

Moreover, for every $\varepsilon>0$, there exists an integer $m$ such that, for
all sufficiently large $t$, there exists a $K_{m+1}$-free graph with
complementary $2$-rank exactly $t$ and at least
\[
t^{\frac{\log 16}{\log 10}-\varepsilon}
\]
vertices.
\end{Proposition}

\begin{proof}
Let $\Gamma$ be the Cayley graph constructed above for $n=5$. Then $\Gamma$ has
$32$ vertices and $144$ edges, and its spectrum is
\[
9^{(1)},\; 5^{(1)},\; 1^{(21)},\; (-3)^{(7)},\; (-7)^{(2)}.
\]
The multiplicity of the eigenvalue $1$ gives
\[
\operatorname{rank}_{\mathbb R}(A_\Gamma-I)=11.
\]
Since the binary rank of an integral matrix is at most its rank over
$\mathbb R$, we get
\[
\operatorname{rank}_{\mathbb F_2}(A_\Gamma+I)
=
\operatorname{rank}_{\mathbb F_2}(A_\Gamma-I)
\leq 11.
\]
In fact, a direct computation gives
\[
\operatorname{rank}_{\mathbb F_2}(A_\Gamma+I)=10.
\]

We now apply Theorem~\ref{thm:BLS} to this graph with
\[
a_1=32,\qquad a_2=144.
\]
For each fixed integer $m\geq 2$ and each positive integer $h$, the theorem
produces a $K_{m+1}$-free graph $\Delta_h$ with complementary $2$-rank at most
$10^h$ and at least
\[
\frac{1}{m+1}
\left(
\frac{32}{\left(32+144\cdot 2^{m+1}\right)^{1/(m+1)}}
\right)^h
\]
vertices.

Since
\[
10^h = T
\]
implies
\[
\left(
\frac{32}{\left(32+144\cdot 2^{m+1}\right)^{1/(m+1)}}
\right)^h
=
T^{\beta_m},
\]
we obtain, for infinitely many values $T=10^h$, a $K_{m+1}$-free graph with
complementary $2$-rank at most $T$ and at least
\[
T^{\beta_m-o(1)}
\]
vertices.

It remains only to pass from the sequence $T=10^h$ to every sufficiently large
$t$. Given $t$, choose $h$ such that
\[
10^h\leq t<10^{h+1}.
\]
Take the graph $\Delta_h$ and add isolated vertices until the complementary
$2$-rank becomes exactly $t$. Adding an isolated vertex adds a $1\times 1$
identity block to $A+I$, so it increases the complementary $2$-rank by exactly
one and does not create any clique. Hence the resulting graph is still
$K_{m+1}$-free, has complementary $2$-rank exactly $t$, and has at least
\[
(10^h)^{\beta_m-o(1)}
\geq
(t/10)^{\beta_m-o(1)}
=
t^{\beta_m-o(1)}
\]
vertices.

Finally,
\[
\lim_{m\to\infty}
\frac{32}{\left(32+144\cdot 2^{m+1}\right)^{1/(m+1)}}
=
16.
\]
Therefore
\[
\lim_{m\to\infty}\beta_m
=
\frac{\log 16}{\log 10}.
\]
Thus, for every $\varepsilon>0$, one can choose $m$ large enough so that
\[
\beta_m>
\frac{\log 16}{\log 10}-\varepsilon,
\]
which proves the final claim.
\end{proof}
This implies Corollary \ref{cor:partovoidW}.

\begin{Remark}
    In \cite{Chawin24}, it is shown that there is an $m$-nearly orthogonal set in $\mathbb{F}_2^t$ of size $t^{\Omega(m/\log m)}$, for $m \rightarrow \infty$. 
    While asymptotically this is a better construction than ours, for small $m$ the constant in $\Omega(m/\log m)$ makes it worse for values like $m = 28$. 
\end{Remark}

\section{Conclusion}

We have established a new connection between partial $m$-ovoids in finite classical polar spaces and Ramsey numbers (Theorem~\ref{Ramseypartialmovoid}). 
It is conceivable that this connection can lead to new explicit constructions for various Ramsey numbers. 
Of particular interest is the case when $m$ is linear in the rank of a polar space, as it can be used to show that $R(t, t) > \mathrm{exp}(t)$ via an explicit construction (see Theorem~\ref{probconstruction} which shows the existence of these objects), which is one of the main open problems in Ramsey theory.
Our connection to partial $m$-ovoids can also lead to improved lower bounds on some small Ramsey numbers. 

Our constructions in the binary symplectic space use graphs with certain rank properties to construct large partial $m$-ovoids. 
It will be of interest to look for similar constructions in other finite classical polar spaces, and in other finite geometrical problems.
For other fields, while we lose the direct connection between partial $m$-ovoids and $m$-nearly orthogonal sets, it might still be possible to use the construction in \cite{Haviv} over arbitrary fields.
We also believe that the probabilistic construction can be improved further. 

Finally, it would be very interesting to improve the upper bounds on partial $m$-ovoids in the binary symplectic space and to give better constructions. 
We also propose determining whether there exists a triangle-free graph $G$ on $N$ vertices with $\mathrm{rank}(A_G + I) = o(N)$.

\subsection*{Acknowledgements}

The third author was supported by a 
postdoctoral fellowship of the Research Foundation -- Flanders (FWO).
The fourth author was supported by an NWO open competition grant (OCENW.M.22.090). 
The first author thanks Delft University of Technology for their hospitality while he was
on sabbatical there, where most of this work was done.
We thank Peter Frankl for discussing the Oddtown Theorem
and its possible generalisations with us.
We also thank Sam Mattheus for various fruitful discussions.

\end{document}